\documentclass{amsart}
\usepackage{latexsym}
\usepackage{amssymb}
\usepackage{amsmath}

\begin{document}


\newtheorem{theorem}{Theorem}
\newtheorem{problem}{Problem}
\newtheorem{definition}{Definition}
\newtheorem{lemma}{Lemma}
\newtheorem{proposition}{Proposition}
\newtheorem{corollary}{Corollary}
\newtheorem{example}{Example}
\newtheorem{conjecture}{Conjecture}
\newtheorem{algorithm}{Algorithm}
\newtheorem{exercise}{Exercise}
\newtheorem{xample}{Example}
\newtheorem{remarkk}{Remark}

\newcommand{\be}{\begin{equation}}
\newcommand{\ee}{\end{equation}}
\newcommand{\bea}{\begin{eqnarray}}
\newcommand{\eea}{\end{eqnarray}}
\newcommand{\beq}[1]{\begin{equation}\label{#1}}
\newcommand{\eeq}{\end{equation}}
\newcommand{\beqn}[1]{\begin{eqnarray}\label{#1}}
\newcommand{\eeqn}{\end{eqnarray}}
\newcommand{\beaa}{\begin{eqnarray*}}
\newcommand{\eeaa}{\end{eqnarray*}}
\newcommand{\req}[1]{(\ref{#1})}

\newcommand{\lip}{\langle}
\newcommand{\rip}{\rangle}
\newcommand{\uu}{\underline}
\newcommand{\oo}{\overline}
\newcommand{\La}{\Lambda}
\newcommand{\la}{\lambda}
\newcommand{\eps}{\varepsilon}
\newcommand{\om}{\omega}
\newcommand{\Om}{\Omega}
\newcommand{\ga}{\gamma}
\newcommand{\ka}{\kappa}
\newcommand{\rrr}{{\Bigr)}}
\newcommand{\qqq}{{\Bigl\|}}

\newcommand{\dint}{\displaystyle\int}
\newcommand{\dsum}{\displaystyle\sum}
\newcommand{\dfr}{\displaystyle\frac}
\newcommand{\bige}{\mbox{\Large\it e}}
\newcommand{\integers}{{\Bbb Z}}
\newcommand{\rationals}{{\Bbb Q}}
\newcommand{\reals}{{\rm I\!R}}
\newcommand{\realsd}{\reals^d}
\newcommand{\realsn}{\reals^n}
\newcommand{\NN}{{\rm I\!N}}
\newcommand{\DD}{{\rm I\!D}}
\newcommand{\degree}{{\scriptscriptstyle \circ }}
\newcommand{\dfn}{\stackrel{\triangle}{=}}
\def\complex{\mathop{\raise .45ex\hbox{${\bf\scriptstyle{|}}$}
     \kern -0.40em {\rm \textstyle{C}}}\nolimits}
\def\hilbert{\mathop{\raise .21ex\hbox{$\bigcirc$}}\kern -1.005em {\rm\textstyle{H}}} 
\newcommand{\RAISE}{{\:\raisebox{.6ex}{$\scriptstyle{>}$}\raisebox{-.3ex}
           {$\scriptstyle{\!\!\!\!\!<}\:$}}} 

\newcommand{\hh}{{\:\raisebox{1.8ex}{$\scriptstyle{\degree}$}\raisebox{.0ex}
           {$\textstyle{\!\!\!\! H}$}}}

\newcommand{\OO}{\won}
\newcommand{\calA}{{\mathcal A}}
\newcommand{\calB}{{\mathcal B}}
\newcommand{\calC}{{\cal C}}
\newcommand{\calD}{{\cal D}}
\newcommand{\calE}{{\cal E}}
\newcommand{\calF}{{\mathcal F}}
\newcommand{\calG}{{\cal G}}
\newcommand{\calH}{{\cal H}}
\newcommand{\calK}{{\cal K}}
\newcommand{\calL}{{\mathcal L}}
\newcommand{\calM}{{\mathcal M}}
\newcommand{\calO}{{\cal O}}
\newcommand{\calP}{{\cal P}}
\newcommand{\calT}{{\mathcal T}} 
\newcommand{\calU}{{\mathcal U}}
\newcommand{\calX}{{\cal X}}
\newcommand{\calY}{{\mathcal Y}}
\newcommand{\calZ}{{\mathcal Z}}
\newcommand{\calXX}{{\cal X\mbox{\raisebox{.3ex}{$\!\!\!\!\!-$}}}}
\newcommand{\calXXX}{{\cal X\!\!\!\!\!-}}
\newcommand{\gi}{{\raisebox{.0ex}{$\scriptscriptstyle{\cal X}$}
\raisebox{.1ex} {$\scriptstyle{\!\!\!\!-}\:$}}}
\newcommand{\intsim}{\int_0^1\!\!\!\!\!\!\!\!\!\sim}
\newcommand{\intsimt}{\int_0^t\!\!\!\!\!\!\!\!\!\sim}
\newcommand{\pp}{{\partial}}
\newcommand{\al}{{\alpha}}
\newcommand{\sB}{{\cal B}}
\newcommand{\sL}{{\cal L}}
\newcommand{\sF}{{\cal F}}
\newcommand{\sE}{{\cal E}}
\newcommand{\sX}{{\cal X}}
\newcommand{\R}{{\rm I\!R}}
\renewcommand{\L}{{\rm I\!L}}
\newcommand{\vp}{\varphi}
\newcommand{\N}{{\rm I\!N}}
\def\ooo{\lip}
\def\ccc{\rip}
\newcommand{\ot}{\hat\otimes}
\newcommand{\rP}{{\Bbb P}}
\newcommand{\bfcdot}{{\mbox{\boldmath$\cdot$}}}

\renewcommand{\varrho}{{\ell}}
\newcommand{\dett}{{\textstyle{\det_2}}}
\newcommand{\sign}{{\mbox{\rm sign}}}
\newcommand{\TE}{{\rm TE}}
\newcommand{\TA}{{\rm TA}}
\newcommand{\E}{{\rm E\,}}
\newcommand{\won}{{\mbox{\bf 1}}}
\newcommand{\Lebn}{{\rm Leb}_n}
\newcommand{\Prob}{{\rm Prob\,}}
\newcommand{\sinc}{{\rm sinc\,}}
\newcommand{\ctg}{{\rm ctg\,}}
\newcommand{\loc}{{\rm loc}}
\newcommand{\trace}{{\,\,\rm trace\,\,}}
\newcommand{\Dom}{{\rm Dom}}
\newcommand{\ifff}{\mbox{\ if and only if\ }}
\newcommand{\nproof}{\noindent {\bf Proof:\ }}
\newcommand{\remark}{\noindent {\bf Remark:\ }}
\newcommand{\remarks}{\noindent {\bf Remarks:\ }}
\newcommand{\note}{\noindent {\bf Note:\ }}

\newcommand{\boldx}{{\bf x}}
\newcommand{\boldX}{{\bf X}}
\newcommand{\boldy}{{\bf y}}
\newcommand{\boldR}{{\bf R}}
\newcommand{\uux}{\uu{x}}
\newcommand{\uuY}{\uu{Y}}

\newcommand{\limn}{\lim_{n \rightarrow \infty}}
\newcommand{\limN}{\lim_{N \rightarrow \infty}}
\newcommand{\limr}{\lim_{r \rightarrow \infty}}
\newcommand{\limd}{\lim_{\delta \rightarrow \infty}}
\newcommand{\limM}{\lim_{M \rightarrow \infty}}
\newcommand{\limsupn}{\limsup_{n \rightarrow \infty}}

\newcommand{\ra}{ \rightarrow }

\newcommand{\ARROW}[1]
  {\begin{array}[t]{c}  \longrightarrow \\[-0.2cm] \textstyle{#1} \end{array} }

\newcommand{\AR}
 {\begin{array}[t]{c}
  \longrightarrow \\[-0.3cm]
  \scriptstyle {n\rightarrow \infty}
  \end{array}}

\newcommand{\pile}[2]
  {\left( \begin{array}{c}  {#1}\\[-0.2cm] {#2} \end{array} \right) }

\newcommand{\floor}[1]{\left\lfloor #1 \right\rfloor}

\newcommand{\mmbox}[1]{\mbox{\scriptsize{#1}}}

\newcommand{\ffrac}[2]
  {\left( \frac{#1}{#2} \right)}

\newcommand{\one}{\frac{1}{n}\:}
\newcommand{\half}{\frac{1}{2}\:}

\def\le{\leq}
\def\ge{\geq}
\def\lt{<}
\def\gt{>}

\def\squarebox#1{\hbox to #1{\hfill\vbox to #1{\vfill}}}
\newcommand{\nqed}{\hspace*{\fill}
          \vbox{\hrule\hbox{\vrule\squarebox{.667em}\vrule}\hrule}\bigskip}

\title{Entropic Solution of the Innovation Conjecture of T. Kailath}

\author{Ali  S\"uleyman  \"Ust\"unel}

\begin{abstract}
\noindent
On a general filtered probability space $(\Om,\calF,(\calF_t,t\in [0,1]),P)$, for a given signal
$U_t=B_t+\int_0^t\dot{u}_sds$, where $(B_t,t\in [0,1])$ is a Brownian
motion and $\dot{u}$ is adapted and in $L^2(dt\times dP)$, we prove that the filtration of $U$,
noted $(\calU_t,t\in [0,1])$,  is
equal to the filtration of its innovation process $Z$ which defined as
$Z_t=U_t-\int_0^t E_P[\dot{u}_s|\calU_s]ds,\,t\in[0,1]$, 
 if and only if 
$$
H(Z(\nu)|\mu)=\half E_\nu\left[\int_0^1|E_P[\dot{u}_s|\calU_s]|^2ds\right]
$$
where $d\nu=\exp(-\int_0^1 E_P[\dot{u}_s|\calU_s]dZ_s-\half
\int_0^1|E_P[\dot{u}_s|\calU_s]|^2 ds)dP$ in case the density has
expectation one, otherwise we give a localized version of the same
strength with a sequence of stopping times of the filtration of $U$.

\end{abstract}
\maketitle
\section{\bf{Introduction}}
\noindent
Let $(\Omega,\calF,(\calF_t,t\in [0,1]),P)$ be a probability space
satisfying the usual conditions
and denote by   $(W,H,\mu)$   the classical Wiener space, i.e.,
$W=C_0([0,1],\R^d)$, $H$ is the corresponding Cameron-Martin space
consisting of $\R^d$-valued  absolutely continuous functions
on $[0,1]$ with square integrable derivatives, which is a Hilbert
space under the norm $|h|_H^2=\int_0^1|\dot{h}(s)|^2ds$, where
$\dot{h}$ denotes the Radon-Nikodym derivative of the absolutely
continuous function $t\to h(t)$ w.r.t. the Lebesgue measure on $[0,1]$. Denote by  $(\calB_t,\,t\in
[0,1])$  the filtration of the canonical Wiener process, completed
w.r.t. $\mu$-negligeable sets. The question that we address in this
paper is the following: assume that $U: \Omega\to W$ is a map of the
following form:
$$
U(\om)(t)=U_t(\om)=B_t(\om)+\int_0^t\dot{u}_s(\om)ds\,,
$$
where $B=(B_t,t\in[0,1])$ is a Brownian motion on $\Om$,
$(s,\om)\to\dot{u}_s(\om)$ is an $\R^d$-valued map belonging to the
space  $L^2_a(P;H)$ which consists of the elements of
$L^2([0,1]\times\Om,\calB([0,1])\otimes \calF,dt\times dP)$ which are
$(\calF_s,s\in[0,1])$-adapted for almost all $s\in [0,1]$.  Let us
define the innovation process $Z$ associated to $U$ as to be
$$
Z_t=U_t-\int_0^t E_P[\dot{u}_s|\calU_s]ds\,,
$$
where $(\calU_t,t\in[0,1])$ is the filtration generated by $U$.  It is
well-known that $Z$ is a $P$-Brownian motion w.r.t. $(\calU_t,t\in[0,1])$.  $Z$
is naturally adapted to $(\calU_t,\,t\in [0,1])$.  This means that the information obtained via
$Z$ is included in the information obtained from $U$. P. Frost \cite{Fr} and  T. Kailath \cite{Kailath} have
conjectured that  in practical situations the converse of this
observation is also true. In \cite{Benes}, V.A. Bene\v{s} has remarked
that this conjecture holds if and only if there is a hidden process
which is a strong solution of a certain stochastic differential
equation from which one can construct the initial system.
This conjecture has also  been proved under restrictive supplementary
hypothesis (cf. \cite{A-M}) where $\dot{u}$ is independent of the
Brownian motion $B$. The main objection to these works lies in the
fact that the condition of \cite{Benes} is unverifiable from the
observed data, hence numerically it is not useful, the second one uses
a hypothesis of independence  which is too strong to be encountered in
the engineering applications. 
In this paper we give a necessary and sufficient condition in the most
general case using the entropic characterization 
of the almost sure invertibility of adapted perturbations of identity
(API in short). Let us explain it briefly for the reader to understand
the idea and the difference from the other works: for simplicitiy,
assume that 
\begin{equation}
\label{integ}
E_P\left[\exp\left(-\int_0^1E_P[\dot{u}_s|\calU_s]dZ_s-\half\int_0^1
    |E_P[\dot{u}_s|\calU_s]|^2ds\right)\right]=1\,,
\end{equation}
and denote by $\rho(-\delta_Z\hat{u})$ the Girsanov exponential inside
the above expectation, here we use the notation $\delta_Z$ to denote  the
It\^o integral of the Lebesgue density of the vector field which is
defined as
$(t,\om)\to\hat{u}(t,\om)=\int_0^tE_P[\dot{u}_s|\calU_s](\om)ds$,  w.r.t. $Z$.
Define a new measure $\nu$ by
$d\nu=\rho(-\delta_Z\hat{u}) dP$. Then the observation process $U$ is
adapted to the filtration of the innovation process $Z$ up to
negligeable sets if and only if we have 
$$
H(Z(\nu)|\mu)=\half E_\nu[|\hat{u}|_H^2]\left(=\half
E_\nu\int_0^1|E_P[\dot{u}_s|\calU_s]|^2ds\right)\,,
$$
where $Z(\nu)$ denotes the push forward of the measure $\nu$ under $Z$
and $H(Z(\nu)|\mu)$ is the relative entropy of $Z(\nu)$ w.r.t. the
Wiener measure $\mu$, i.e., 
$$
H(Z(\nu)|\mu)=\int_W \frac{dZ(\nu)}{d\mu}\log \frac{dZ(\nu)}{d\mu}
d\mu\,.
$$
As it is clear, the verification of  this condition, namely  the equality of the entropy to the
total kinetic energy of $\hat{u}$, requires only the knowledge about
the observation process $U$. However the calculation of the relative
entropy may be time consuming. In fact, as it follows from Theorem
\ref{stopping-thm}, all these results are valid when one works causally
with time, in other words, they hold also when one works on the time
interval $[0,t]$, for $t\geq 0$ since they are restictable even to the
random time intervals. The final result says that we can also suppress
the hypothesis (\ref{integ}) using a sequence of $(\calU_t,t\in
[0,1])$-stopping times.

\section{\bf{Characterization of the invertible shifts on the canonical space}}
\noindent
We  begin with the definition of the notion of  almost sure
invertibility with respect to a measure. This notion is extremely
important since it makes the  things work. Let us note that in this
section all the expectations and conditional expectations are taken
w.r.t. the Wiener measure $\mu$.

\begin{definition}
 Let $T:W\to W$ be a measurable map.
\begin{itemize}
\item T is called ($\mu$-) almost surely left invertible
if there exists a measurable  map $S:W\to W$ such that  $S\circ T =I_W$ $\mu$-a.s.
\item  Moreover, in this case it is trivial to
see that $T\circ S=I_W$ $T\mu$-a.s., where $T\mu$ denotes the image of
the measure $\mu$ under the map $T$. 
\item If $T\mu$ is equivalent to $\mu$,
then we say in short that $T$ is $\mu$-a.s. invertible. 
\item Otherwise, we
may say that $T$ is $(\mu,T\mu)$-invertible in case  precision is
required or just $\mu$-a.s. left invertible and $S$ is called the
$\mu$-left inverse of $T$.
\end{itemize}
\end{definition}
\noindent
Let $L^2_a(\mu,H)$ be the $\mu$-square integrable  equivalence classes of Cameron-Martin
space (denoted by $H$)-valued functions, hence  $t\to(w\to u(w))$
is an absolutely continuous  function of $t\in [0,1]$ with a
$ds$-square integrable Lebesgue
density denoted as $\dot{u}_s(w)$, moreover we assume that
$w\to\dot{u}_s(w)$ is $\calB_s$-measurable for $ds$-almost all $s\in
[0,1]$, which is a Hilbert space. For short we call them adapted vector fields of class $L^2$.
Similarly $L^0_a(\mu,H)$ denotes the set of adapted vector fields
whose Cameron-Martin norm is $\mu$-a.s. finite, under the topology of
convergence in probability, $L^0_a(\mu,H)$ is a non- locally convex
Fr\'echet space. For reader's convenience, let us note
that $L^0_a(\mu,H)$ is the completion of $L^2_a(\mu,H)$ under the
topology of convergence in probability.

Although the following theorem has been proved in \cite{ASU-2}, for the
reader's convenience we give a short and different proof:
\begin{theorem}
\label{thm-0}
For any $u\in L_a^2(\mu,H)$, we have the following inequality
$$
H(U\mu|\mu)\leq \frac{1}{2}E\int_0^1|\dot{u}_s|^2ds\,,
$$
where $H(U\mu|\mu)$ is the relative entropy of the measure $U\mu$
w.r.t. $\mu$. 
\end{theorem}
\nproof
Let $L$ be the Radon-Nikodym density of $U\mu$ w.r.t. $\mu$. For any
$0\leq g\in C_b(W)$, using the Girsanov theorem, we have 
$$
E[g\circ U]=E[g\,L]\geq E[g\circ U\,L\circ U\,\rho(-\delta u)]\,,
$$
hence 
$$
L\circ U\,E[\rho(-\delta u)|U]\leq 1
$$
$\mu$-a.s. Consequently, using the Jensen inequality
\beaa
H(U\mu|\mu)&=&E[L\log L]=E[\log L\circ U]\\
&\leq&-E[\log E[\rho(-\delta u)|U]]\\
&\leq&-E[\log \rho(-\delta u)]\\
&=&\half E\int_0^1|\dot{u}_s|^2ds\,.
\eeaa
\nqed

\begin{theorem}
\label{main-thm}
Assume that $U=I_W+u$ is an API, i.e., $u\in L^2_a(\mu,H)$ such that
$s\to\dot{u}(s,w)$ is $\calB_s$-measurable for almost all $s$. 
Then $U$ is almost surely left invertible with a left inverse $V$ if and only if 
$$
H(U\mu|\mu)=\frac{1}{2}E[|u|_H^2]=\frac{1}{2}E\int_0^1|\dot{u}_s|^2ds\,,
$$
i.e., if and only if the entropy of $U\mu$ is equal to the energy of
the drift $u$.
\end{theorem}
\nproof
Due to Theorem \ref{thm-0}, the relative entropy is finite as soon as 
$u\in L^2_a(\mu,H)$. Let us suppose now that the equality holds and let us denote by $L$ the Radon-Nikodym derivative of
$U\mu$ w.r.t. $\mu$. Using the It\^o representation theorem, we can
write 
$$
L=\exp\left(-\int_0^1 \dot{v}_sdW_s-\half\int_0^1|\dot{v}_s|^2ds\right)
$$
$U\mu$-almost surely. Let $V=I_W+v$, as described in \cite{Foll}, from the It\^o formula and Paul
L\'evy's theorem, it is immediate that $V$ is an $U\mu$-Wiener
process, hence 
\begin{equation}
\label{ent}
E[L\log L]=\half E[L\,|v|_H^2]\,.
\end{equation}
Now, for any $f\in C_b(W)$, we have from the Girsanov
theorem
$$
E[f\circ U]=E[f\,L]\geq E[f\circ U\,L\circ U\,\rho(-\delta u)]
$$
consequently
$$
L\circ U\,E[\rho(-\delta u)|U]\leq 1
$$
$\mu$-a.s. Let us denote $E[\rho(-\delta u)|U]$ by $\hat{\rho}$. We
have then $\log L\circ U+\log\hat{\rho}\leq 0$ $\mu$-a.s. Taking the
expectation w.r.t. $\mu$ and the Jensen inequality give
\beaa
H(U\mu|\mu)&=&E[L\log L]\leq -E[\log \hat{\rho}]\\
&\leq&-E[\log\rho(-\delta u)]=\half E[|u|_H^2]\,.
\eeaa
Since $\log$ is a strictly concave function, the equality $E[\log
\hat{\rho}]=E[\log\rho(-\delta u)]$ implies that $\rho(-\delta u)=\hat{\rho}$
$\mu$-a.s. Hence we obtain
$$
E[L\log L+\log\rho(-\delta u)]=E[\log(L\circ U\,\rho(-\delta u))]=0\,,
$$
and since $L\circ U\rho(-\delta u)\leq 1$ $\mu$-a.s., we should have 
\begin{equation}
\label{basic}
L\circ U\rho(-\delta u)=1
\end{equation}
$\mu$-a.s. Combinrning  the exponential representation of $L$ with the
relation (\ref{basic}) implies
\bea
\label{basic-1}
0&=&\left(\int_0^1\dot{v}_sdW_s\right)\circ U+\half|v\circ U|_H^2+\delta
u+\half|u|_H^2\nonumber\\
&=&\delta(v\circ U)+\delta u+(v\circ U,u)_H+\half(|u|_H^2+|v\circ
U|_H^2)\nonumber\\
&=&\delta(v\circ U+u)+\half|v\circ U+u|_H^2
\eea
$\mu$-a.s. From the relation (\ref{ent}) it follows that $v\circ U\in
L^2_a(\mu,H)$, hence taking the expectations on both sides of
(\ref{basic-1}) w.r.t. $\mu$ is licit  and this implies $v\circ U+u=0$ $\mu$-a.s., which means that
$V=I_W+v$ is the $\mu$-left inverse of $U$. 

To show the neccessity, let us denote by $(L_t,t\in [0,1])$ the
martingale
$$
L_t=E[L|\calB_t]=E\left[\frac{dU\mu}{d\mu}|\calB_t\right]
$$
and let 
$$
T_n=\inf\left(t:\,L_t<\frac{1}{n}\right)\,.
$$
Since $U\circ V=I_W$ $(U\mu)$-a.s.,  $V$ can be written as $V=I_W+v$
($U\mu$)-a.s. and that $v\in L_a^0(U\mu,H)$, i.e., 
$v(t,w)=\int_0^t \dot{v}_s(w)ds$, $\dot{v}$
is adapted to the filtration $(\calB_t)$ completed w.r.t. $U\mu$ and
$\int_0^1|\dot{v}_s|^2ds<\infty$ $(U\mu)$-a.s. Noting that $\{t\leq
T_n\}\subset \{L_t>0\}$ and hence that $\mu$ and $U\mu$
are equivalent on $\calF_{T_n}$, we conclude
$$
\int_0^{T_n}|\dot{v}_s|^2ds<\infty
$$
{\bf{$\mu$-almost surely}}. Consequently the  inequality 
$$
E_\mu[\rho(-\delta v^n)]\leq 1
$$
holds true for any $n\geq 1$, where $v^n(t,w)=\int_0^t
1_{[0,T_n]}(s,w)\dot{v}_s(w)ds$. By positivity we also have 
$$
E_\mu[\rho(-\delta v^n)1_{\{L>0\}}]\leq 1\,.
$$
Since $\lim_{n\to\infty}T_n=\infty $ $(U\mu)$-a.s., we also have $\lim_{n\to\infty}T_n=\infty
$ $\mu$-a.s. on the set $\{L>0\}$ and the Fatou lemma implies
\begin{equation}
\label{ineq-1}
E_\mu[\rho(-\delta v)1_{\{L>0\}}] =E_\mu[\lim_n\rho(-\delta
v^n)1_{\{L>0\}}]\leq \lim\inf_n E_\mu[\rho(-\delta v^n)1_{\{L>0\}}]\leq 1\,.
\end{equation}
From the identity $U\circ V=I_W$
$(U\mu)$-a.s., we have $v+u\circ V=0$ $(U\mu)$-a.s., hence $v\circ
U+u=0$ $\mu$-a.s.  An algebraic  calculation gives immediately 
\begin{equation}
\label{ineq-2}
\rho(-\delta v)\circ U\,\rho(-\delta u)=1
\end{equation}
$\mu$-a.s. Now applying the Girsanov theorem to API $U$ and using the
relation (\ref{ineq-2}), we obtain 
\beaa
E[g\circ U]&=&E[g\,L]=E\left[g\circ U(\rho(-\delta v)1_{\{L>0\}})\circ
U\rho(-\delta u)\right]\\
&\leq&E\left[g\,\rho(-\delta v)1_{\{L>0\}}\right]\,,
\eeaa
for any positive $g\in C_b(W)$ (note that on the set $\{L>0\}$, 
$\rho(-\delta v)$ is perfectly well-defined w.r.t. $\mu$). Therefore 
$$
L\leq \rho(-\delta v)1_{\{L>0\}}
$$
$\mu$-a.s. Now, this last inequality, combined with the inequality
(\ref{ineq-1}) entails that 
$$
L=\rho(-\delta v)1_{\{L>0\}}
$$
$\mu$-a.s., hence 
$$
L\circ U\,\rho(-\delta u)=1
$$
$\mu$-a.s. To complete the proof it suffices to remark then that 
\beaa
H(U\mu|\mu)&=&E[L\log L]=E[\log L\circ U]\\
&=&E[-\log \rho(-\delta u)]\\
&=&
\half E[|u|_H^2]\,.
\eeaa
\nqed

\noindent
The following result comes almost for free:
\begin{theorem}
\label{stopping-thm}
Assume that $U=I_W+u$ is an API which is $\mu$-a.s. left invertible,
let $\tau$ be any stopping time such that $u^\tau$, defined as
$u^\tau(t,w)=u(t\wedge\tau(w),w)$ such that 
$$
E[\rho(-\delta u^\tau)]=1\,.
$$
Then 
$U^\tau=I_W+u^\tau$ is $\mu$-a.s. invertible, in other words there
exists some API, say $V'$ such that $V'\circ U^\tau=U^\tau\circ V'=I_W$
$\mu$-a.s.
\end{theorem}
\nproof
Since $E[\rho(-\delta u^\tau)]=1$, $U^\tau\mu$ is equivalent to the
Wiener measure $\mu$, hence its Radon-Nikodym density can be written as 
$$
\frac{dU^\tau\mu}{d\mu}=\rho(-\delta \xi)\,.
$$
From the Girsanov theorem it follows that 
\begin{equation}
\label{M-A}
\rho(-\delta\xi)\circ U^\tau\,E[\rho(-\delta u^\tau)|U^\tau]=1
\end{equation}
$\mu$-a.s. Let $z$ be the innovation process of $U^\tau$, which is
defined as $z_t=U^\tau_t-\int_0^t E[\dot{u}^\tau_s|\calU^\tau_s]ds$,
where $(\calU^\tau_s,s\in [0,1])$ denotes the filtration corresponding
to $U^\tau$. Applying the
Girsanov theorem again, this time using the Brownian motion $z$
(cf. \cite{ASU-2} for the details), we
find that 
$$
E[\rho(-\delta
u^\tau)|U^\tau]=\exp\left(-\int_0^1E[\dot{u}^\tau_s|\calU^\tau_s]dz_s-\half
\int_0^1|E[\dot{u}^\tau_s|\calU^\tau_s]|^2 ds\right)\,.
$$
This relation, combined with the equation (\ref{M-A}) gives the
relation 
$$
\dot{\xi}_t\circ U^\tau+ E[\dot{u}\,1_{[0,\tau]}(t)|\calU^\tau_t]=0
$$
$dt\times d\mu$-a.s. Besides, for any $A\in L^\infty(\mu)$, we have 
\beaa
E[A\,E[\dot{u}_t\,1_{[0,\tau]}(t)|\calU^\tau_t]]&=&E[E[A|\calU^\tau_t]\dot{u}_t\,1_{[0,\tau]}(t)]\\
&=&E[E[A|\calU_t]\dot{u}_t\,1_{[0,\tau]}(t)]\\
&=&E[A\,E[\dot{u}^\tau_t|\calU_t]]\\
&=&E[A\,\dot{u}^\tau_t]\,,
\eeaa
where the last equality follows from the left invertibility of $U$. 
Hence we obtain
$$
\dot{\xi}_t\circ U^\tau+ \dot{u}\,1_{[0,\tau]}(t)=\dot{\xi}_t\circ U^\tau+ \dot{u}^\tau_t=0
$$
$dt\times d\mu$-a.s., which is equivalent to $\mu$-a.s. invertibility
of $U^\tau$.
\nqed

\section{\bf{The case of a general probability space}}
\noindent
The following result is essential for the proof of the conjecture where 
we use the notations explained in the introduction and we
differentiate carefully the Wiener measure $\mu$  and the probability
$P$ as well as the  respective  expectations and conditional
expectations to avoid any ambiguity.  In particular, we denote by
$L^2_a(P,H)$ the space of adapted, $P$-square integrable vector
fields, this is exactly the same space as $L^2_a(\mu,H)$, where the
Wiener space is replaced by a general probability space $\Om$, $\mu$
is replaced by a probability $P$ defined on $(\Om,\calF)$, and the
canonical filtration $(\calB_t,t\in [0,1])$ of the Wiener space is replaced by a general filtration
$(\calF_t,t\in [0,1])$ of $(\Om,\calF)$. Similarly $L^0_a(P,H)$
denotes the version of $L^0_a(\mu,H)$ on this general probability
space, remark that it is the completion of $L^2_a(P,H)$
w.r.t. convergence in probability $P$.
\begin{theorem}
\label{gen-thm}
Let $U=B+u=B+\int_0^\cdot \dot{u}_sds$ be an adapted perturbation of identity mapping $\Om$ to $W$ with $u\in L^2(P,H)$. Then
$$
H(U(P)|\mu)=\half E_P[|u|_H^2]
$$
if and only if there exists some $v:W\to H$ (of the form $v=\int_0^\cdot \dot{v}_s ds$) with $\dot{v}$ adapted $ds$-a.s. to the filtration $(\calB_t(W))$ such that 
$$
U(\om)=B(\om)-v\circ U(\om)\,,
$$
which implies in particular that $B=Z$, where $Z$ is the innovation process associated to $U$, in other words $U$ is a solution of the following stochastic differential equation 
$$
dU_t=-\dot{v}_t\circ U dt+dB_t\,.
$$
\end{theorem}
\nproof Note first that $U$ is not neccessarily a strong solution. Let us now prove the neccessity: Since $U$ is an API, $U(P)$ is absolutely continuous w.r.t. the Wiener measure $\mu$, let $l$ be the corresponding Radon-Nikodym derivative. We can represent it as a Girsanov exponential $U(P)$-a.s., i.e., we have
\beaa
l&=&\frac{dU(P)}{d\mu}=\rho(-\delta v)\\
&=&\exp\left(-\int_0^1\dot{v}_sdW_s-\half\int_0^1|\dot{v}_s|^2ds\right)\,,
\eeaa
$U(P)$-a.s., where $(W_t)$ is the canonical Wiener process. For any positive  $f\in C_b(W)$, it follows from the Girsanov theorem
$$
E_P[f\circ U]=E_\mu[f\,l]\geq E_P[f\circ U\,l\circ U\,\rho(-\delta_B u)]\,,
$$
where
$$
\rho(-\delta_B u)=\exp\left(-\int_0^1\dot{u}_sdB_s-\half\int_0^1|\dot{u}_s|^2ds\right)\,.
$$
This inequality, which is valid for any positive, measurable $f$, implies that 
$$
l\circ U\,E_P[\rho(-\delta_B u)|U]\leq 1
$$
$P$-a.s. Therefore
\beaa
H(U(P)|\mu)&=&E_\mu[l\log l]=E_P[\log l\circ U]\\
&\leq&-E_P[\log E_P[\rho(-\delta_B u)|U]]\leq \half E_P[|u|_H^2]\,.
\eeaa
The equality hypothesis $H(U(P)|\mu)=\half E_P[|u|_H^2]$ and the strict convexity of the function $x\to-\log x$ imply that 
$$
l\circ U\,\rho(-\delta_B u)=1
$$
$P$-a.s. Therefore
\beaa
1&=&\rho(-\delta v)\circ U\rho(-\delta B)\\
&=&\exp-\left[(\delta v)\circ U+\half|v\circ U|_H^2+\delta_B u+\half|u|_H^2\right]\\
&=&\exp-\left[\delta_B (v\circ U)+(v\circ U,u)_H+\half|v\circ U|_H^2+\delta_B u+\half|u|_H^2\right]\,,
\eeaa
which implies that
$$
\delta_B(u+v\circ U)+\half|v\circ U+u|_H^2=0
$$
$P$-a.s. Since $E_P[|v\circ U|_H^2]=E_\mu[l|v|_H^2]=2E_\mu[l\log l]$, it follows that $v\circ U+u=0$ $P$-a.s. Note that we can write
$$
U=B+u=Z+\hat{u},\,\,\hat{u}=\int_0^\cdot E_P[\dot{u}_s|\calU_s]ds\,,
$$
since $u=-v\circ U$, $\dot{u}$ is adapted to the filtration of $U$, therefore $B=Z$.\\
\noindent
{\bf{Sufficiency:}} If $U=B-v\circ U$, then $Z=B$ and $v\circ U+u=0$. Let $l$ denote again the Radon-Nikodym derivative of $U(P)$ w.r.t. $\mu$, as before we can write $l=\rho(-\delta \xi)$ $U(P)$-a.s., for some $\xi:W\to H$ such that $\xi=\int_0^\cdot \dot{\xi}_sds$, $\int_0^1|\dot{\xi}_s|^2ds<\infty$ $U(P)$-a.s. and $ \dot{\xi}_s$ is $\calB_s(W)$-measurable $ds$-a.s. Using the Girsanov theorem as above, we find that 
$$
l\circ U\,E_P[\rho(-\delta_B u)|U]\leq 1
$$
but the hypothesis implies that $\rho(-\delta_B u)$ is $U$-measurable, it then follows that 
$$
\delta_B(u+\xi\circ U)+\half|\xi\circ U+u|_H^2\leq 0
$$
$P$-a.s. Since $E_P[|\xi\circ U|_H^2]=2H(U(P)|\mu)<\infty$, it follows that $\xi\circ U=v\circ U$ $P$-a.s. Consequently
$$
H(U(P)|\mu)=E_\mu[l\log l]=E_P[\log l\circ U]=-E_P[\log\rho(-\delta_B u)]=\half E_P[|u|_H^2]
$$
and this completes the proof.
\nqed

\noindent
Theorem \ref{gen-thm} says that $U=B+u$ with $u\in L^0_a(P,H)$,  is the weak solution of the SDE 
$$
dU_t=dB_t-\dot{v}_t\circ U dt
$$ 
if and only if we have the equality between the
entropy $H(U(P)|\mu)$ and the total kinetic energy of $u$ w.r.t. 
the probability $P$. A natural question is:  when this solution is
strong? The following theorem gives the answer:
\begin{theorem}
\label{strong-sol}
Assume that $U$ is a weak solution of the SDE
$$
dU_t=dB_t-\dot{v}_t\circ U dt\,,
$$ 
with the hypothesis that $v\circ U\in L_a^0(P,H)$, define the sequence
of stopping times $(t_n, n\geq 1)$ as 
$$
t_n=\inf\left(t:\,\int_0^t|\dot{v}_s|^2ds>n\right)
$$
 let 
$$
\dot{u}_n(t,\om)=-\dot{v}_t\circ U(\om)\,1_{[0,t_n\circ U]}(t)
$$
and let $U^n=B+u_n$ where $u_n(t,\om)=\int_0^t\dot{u}_n(s,\om)ds$.
Define a new probability $Q_n$  by
$dQ_n=\rho(-\delta_B(u_n))\,dP$. Then
$$
H(B(Q_n)|\mu)=\half E_{Q_n}[|u_n|_H^2]\,,
$$
for any $n\geq 1$  if and only if 
$U$ is a strong solution.

\end{theorem}
\nproof
{\bf{Necessity:}} Since under $Q_n$, $U_n$ is a Brownian motion and the hypothesis
combined with Theorem \ref{gen-thm} implies that $v_n\circ U$ is
measurable w.r.t. the filtration of $B$ up to $Q_n$-negligeable sets,
since $Q_n$ is equivalent to $P$, it follows that $v_n\circ U$ is
adapted to the same filtration completed with $P$-negligeable
sets. Since $\lim_{n\to\infty}v_n\circ U=v\circ U$, $U$ is also adapted to the
$P$-completion of the filtration of $B$, hence $U$ is a strong
solution of the above SDE.
\newline
{\bf{Sufficiency:}} If $U$ is a strong solution, then it is of the
form $U=\hat{U}(B)=B-v\circ \hat{U}(B)$ and $\hat{U}:W\to W$ has a
$\mu$-a.s. left inverse $V=I_W+v$. Since $Q_n$ is equivalent to $P$,
we have also $U=\hat{U}(B)$ $Q_n$-a.s. Moreover $B=U_n+v_n\circ U$ and
$v_n\circ U$ is adapted to the filtration of $B$ up to
$Q_n$-negligeable sets for any $n\geq 1$. Due to Theorem \ref{gen-thm} this is
equivalent to the equality 
$$
H(B(Q_n)|\mu)=\half E_{Q_n}[|u_n|_H^2]\,,
$$
for any $n\geq 1$.
\nqed

\section{\bf{Proof of the Innovation Conjecture}}
\noindent
We are now at a position to give the proof of the conjecture. We shall
do it in two steps using the notations explained in the introduction.
The first step is with a supplementary hypothesis to explain clearly
the idea, the second one is in full generality. \\
We have the relation 
$$
U=B+u=Z+\hat{u}\,,
$$
and we shall denote by $(\calZ_t,t\in [0,1])$ the filtration generated by the innovation process $Z$. We use also the notation
$$
\rho(-\delta_Z\hat{u})=\exp\left(-\int_0^1 E_P[\dot{u}_s|\calU_s] dZ_s-\half\int_0^1 |E_P[\dot{u}_s|\calU_s]|^2ds\right)\,.
$$
First we give a proof with a supplementary hypothesis which will be suppressed at the final proof:
\begin{proposition}
\label{pf-1}
Assume that 
$$
E_P[\rho(-\delta_Z\hat{u})]=1\,,
$$
denote then by $\nu$ the probability defined by
$d\nu=\rho(\delta_Z\hat{u}) dP$. Then $\calU_t=\calZ_t$ for any $t\geq
0$ up to negligeable sets and $\hat{u}=v\circ Z$, 
with $v\in L^0(\mu,H)$ with $\dot{v}_s$ being $\calB_s(W)$-measurable $ds$-almost surely, if and only if 
$$
H(Z(\nu)|\mu)=\half E_\nu[|\hat{u}|_H^2]\,.
$$
\end{proposition}
\nproof
By Paul L\'evy's Theorem, $U$ is a Brownian motion under the measure
$\nu$ and $Z=U-\hat{u}$. 
Then Theorem \ref{gen-thm} says that (replacing $B$ by $U$ and $P$ by
$\nu$), $\hat{u}$ is a functional of $Z$ and that 
$s\to E_P[\dot{u}_s|\calU_s]$ is adapted to the filtration $(\calZ_s,s\in [0,1])$ $ds$-a.s. Hence $U$ is $Z$-measurable. Moreover, same theorem implies the existence of some $v\in L^0(\mu,H)$ which is defined as 
$$
\frac{dZ(\nu)}{d\mu}=\rho(\delta v)
$$
such that $\hat{u}=v\circ Z$.
\nqed

\noindent
Now we are ready to give the full proof:
\begin{theorem}
\label{inn-thm}
Let $T_n=\inf(t:\,\int_0^t|E_P[\dot{u}_s|\calU_s]|^2ds>n)$, define
\beaa
\hat{u}_n(t,\om)&=&\hat{u}(t\wedge T_n,\om)\\
U_n&=&Z+\hat{u}_n\,.
\eeaa
Then $\calZ_t=\calU_t$ for any $t\geq 0$ up to negligeable sets and
hence $\hat{u}$ should be of the form $\tilde{u}\circ Z$ with some $\tilde{u}\in L^0_a(\mu,H)$ if and only if we have 
\begin{equation}
\label{entropy}
H(Z(\nu_n)|\mu)=\half E_{\nu_n}[|\hat{u}_n|_H^2]
\end{equation}
for any $n\geq 1$, where $d\nu_n=\rho(-\delta_Z\hat{u}_n) dP$, 
and
$$
\rho(-\delta_Z\hat{u}_n)=\exp\left(-\int_0^{T_n}E_P[\dot{u}_s|\calU_s] dZ_s-\half\int_0^{T_n}|E_P[\dot{u}_s|\calU_s]|^2 ds\right)\,.
$$
\end{theorem}
\nproof
{\bf{Sufficiency:}} Under the measure $\nu_n$, $U_n$ is a Brownian
motion and $Z=U_n-\hat{u}_n$. It follows from Theorem \ref{gen-thm}
that $\hat{u}_n$ is $((\calZ_t),\nu_n)$-adapted if and only if the
relation (\ref{entropy}) holds true. Since $\nu_n$ is equivalent to
$P$, $U_n$ is also $((\calZ_t),P)$-adapted for any $n\geq 1$, since
$U_n\to U$ in $L^0(P,W)$ (i.e., $P$-equivalence classes of $W$-valued
measurable functions under the topology of convergence in probability $P$), $U$ is also $((\calZ_t),P)$-adapted.
\\
{\bf{Neccessity:}} Assume that $U$ is $((\calZ_t),P)$-adapted, then it
is also $((\calZ_t),\nu_n)$-adapted since $\nu_n\sim P$ for any $n\geq
1$. Hence $U_n$ is also $((\calZ_t),\nu_n)$-adapted  for any $n\geq 1$
and this is equivalent to the relation (\ref{entropy}) for any $n\geq
1$.
\nqed

\noindent{\bf{Acknowledgement:}} We are grateful to an anonymous
referee who converted this paper into a readable one.


\vspace{2cm}

{\footnotesize{\bf{
\noindent
A.S. \"Ust\"unel, Telecom-Paristech (formerly ENST),
 Dept. Infres,\\
46, rue Barrault, 75013 Paris, France\\
email: ustunel@telecom-paristech.fr}
}}


\begin{thebibliography}{99}




\bibitem{A-M}
D. Allinger and S.K. Mitter: `` New results on the innovations
problem for nonlinear filtering''.  Stochastics  {\bf 4}, no. 4,
339--348, 1980.

\bibitem{Benes}
V. A. Bene\v{s}: ``On Kailath's Innovation Conjecture Hold''. The Bell
System Tech. Jour., Vol. 55, no. 7, 981-1001, 1976.











\bibitem{Foll}
H. F\"ollmer:
``An entropy approach to the time reversal of diffusion processes''. 
Stochastic Differential Systems Filtering and Control
Proceedings of the IFIP-WG 7/1 Working Conference Marseille-Luminy, p.
156-163, France, March 12–17, 1984.

\bibitem{Fr}
P. Frost: ``Estimation in Continuous-Time Nonlinear
Systems''. Dissertation, Stanford University, Stanford, Calif., June 1968.






\bibitem{Kailath}
T. Kailath: ``Some Extensions of the Innovations Theorems''. B.S.T.J.,
50, p. 1487-1494, 1971.

\bibitem{Kailath-Z}
T. Kailath and M. Zakai: ``Absolute Continuity and Radon-Nikodym of
Certain Measures Relative to Wiener Measure''. Ann. Math. stat. 42,
p. 130-140, 1971.




\bibitem{Tsi}
B.S. Tsirelson: ``An example of stochastic differential equation
having no strong solution''. Theor. Prob. Appl. {\bf 20}, p. 416-418, 1975.




\bibitem{ASU}
A. S. \"Ust\"unel:
{\sl Introduction to Analysis on Wiener Space}.
Lecture Notes in Math. Vol. 1610. Springer, 1995.

\bibitem{ASU-1}
A. S. \"Ust\"unel: {\sl Analysis on Wiener Space and Applications}. 
http://arxiv.org/abs/1003.1649, 2010.




\bibitem{ASU-3}
A. S. \"Ust\"unel:``A necessary and sufficient condition for invertibility of
  adapted perturbations of identity on Wiener space''. Comptes Rendus
  Math\'ematiques, Vol. 346, p.  897-900. {\bf 2008}.

\bibitem{FILT}
A. S. \"Ust\"unel and M. Zakai: ``The construction of filtrations on abstract Wiener space''.
J. Funct. Anal. {\bf  143}  , p. 10--32, 1997.


\bibitem{ASU-2} 
A. S. \"Ust\"unel : ``Entropy, invertibility and variational calculus of adapted
  shifts on Wiener space''. Journal of Functional Analysis, Volume
  257, Issue 11, Pages 3655-3689,
{\bf 2009}.


\bibitem{BOOK}
A. S. \"Ust\"unel and M. Zakai:
{\sl Transformation of Measure on Wiener Space}. Springer Verlag, 1999.


\end{thebibliography}
\end{document}